\documentclass[12pt]{article}
\usepackage[a4paper]{geometry}
\usepackage[italian,english]{babel}

\usepackage{amsmath}
\usepackage{amsfonts}
\usepackage{amstext}
\usepackage{amssymb}
\usepackage{amsthm}
\usepackage{amscd}

\usepackage[pagebackref,draft=false]{hyperref}
\hypersetup{colorlinks,
linkcolor=myrefcolor,
citecolor=mycitecolor,
urlcolor=myurlcolor}

\usepackage[capitalize]{cleveref}
\usepackage{caption}

\usepackage{xcolor}
\definecolor{myurlcolor}{rgb}{0,0,0.4}
\definecolor{mycitecolor}{rgb}{0,0.5,0}
\definecolor{myrefcolor}{rgb}{0.5,0,0}
\usepackage{graphicx}
\usepackage{tikz}
\usepackage{tikz-cd}

\usepackage{etoolbox}
\usepackage{makeidx}
\usepackage{authblk} 
\usepackage{sectsty}
\usepackage{mathrsfs}
\usepackage{dsfont}
\usepackage{enumitem} 
\usepackage[]{latexsym}
\usepackage{braket}
\usepackage{caption}
\usepackage[utf8]{inputenx}
\usepackage[T1]{fontenc}
\usepackage{lmodern}
\usepackage{textcomp}
\usepackage{microtype}
\usepackage{totcount}
\usepackage{blindtext}

\newtheorem{example}{Example}

\newtheorem*{proof*}{Proof}

\newcommand{\be}{\begin{equation}}
\newcommand{\ee}{\end{equation}}
\newcommand{\bea}{\begin{eqnarray}}
\newcommand{\eea}{\end{eqnarray}}

\numberwithin{equation}{section}

\newcommand{\M}{\mathcal{M}}
\newcommand{\x}{\mathbf{x}}
\newcommand{\y}{\mathbf{y}}

\renewcommand{\i}{\mathrm{in}}
\newcommand{\f}{\mathrm{fin}}

\title{Generalised Potential Functions in Differential Geometry and Information Geometry}

\author[1,2]{F.M.~Ciaglia}
\author[1,2]{G.~Marmo}
\author[3,4]{J.M.~Pérez-Pardo}

\affil[1]{\textit{\footnotesize Dipartimento di Fisica ``E. Pancini'', Universit\`a di Napoli Federico II, Complesso Universitario di Monte S. Angelo Edificio 6, via Cintia, 80126 Napoli, Italy.}}

\affil[2]{\textit{\footnotesize INFN-Sezione di Napoli, Complesso Universitario di Monte S. Angelo Edificio 6, via Cintia, 80126 Napoli, Italy.}}

\affil[3]{\textit{\footnotesize Depto. de Matem\'aticas, Univ. Carlos III de Madrid, Avda. de la Universidad 30, 28911 Legan\'es, Madrid, Spain.}}

\affil[3]{\textit{\footnotesize Instituto de Ciencias Matem\'aticas (CSIC-UAM-UC3M-UCM), C./ Nicol\'as Cabrera 13-15. Campus de Cantoblanco, UAM, 28049 Madrid, Spain.}}

\date{{\footnotesize Subject Classification 53C99, 53Z05, 53A40} \\ \vskip 1ex {\footnotesize Keywords: \emph{Information Geometry, Generalised Potential Functions, Inverse Problem}}}

\begin{document}


\maketitle

\begin{abstract}
Potential functions can be used for generating potentials of relevant geometric structures for a Riemannian manifold such as the Riemannian metric and affine connections. We study wether this procedure can also be applied to tensors of rank four and find a negative answer. We study this from the perspective of solving the inverse problem and also from an intrinsic point of view.
\end{abstract}



\section{Introduction}

In 1925 R.A.~Fisher \cite{Fisher25} introduced the celebrated concept of statistical distance. 
This notion was the cornerstone for 
what has developed nowadays into a major field in Statistics, known as Information Geometry. 
This has ramifications, for instance, 
in Information Theory, Statistical Mechanics, Biology and Quantum Mechanics, see \cite{AN00, Centcov82}. 
One key ingredient of 
the successful approach of Information Geometry has been the geometrisation. 
That is, aiming  at stating the results in a way that 
is independent of the specific situation or coordinates and thus providing a very general framework that can be applied to a variety 
of situations. 

The central object of the approach is known as a statistical manifold and we will denote it by $\mathcal{S}$. 
This assumes that, given a particular family of probability distributions $\mathcal{P}(X)$ on a measure space $X$, one is able to find a parametrisation, i.e.\ a mapping $$p: \Xi \to \mathcal{P}(X),$$
such that $p$ is a bijection onto its image $p(\Xi)$. Under suitable conditions, we refer to \cite{AN00} for a proper definition of a Statistical Manifold, $
\mathcal{S}$ has the structure of a differentiable manifold and Fisher's statistical distance provides it with a natural Riemannian 
metric, the so called Fisher-Rao metric \cite{Rao45}.

In addition to this Riemannian structure, there are two geometric objects related to it that play an important role. On one hand 
there 
is  a so called dualistic structure \cite{Aea87}. This is a pair of torsion-free connections $(\nabla,\nabla^*)$ that satisfy the relation 
$$X(g(Y,Z)) = g(\nabla_{X}Y, Z) + g(Y,\nabla^*_{X}Z),$$
where $X, Y, Z$ are vector fields over the manifold.
These two connections are said to be dual to each other. Notice that, in particular, the Levi-Civita connection is the unique self--dual connection.
On the other hand there is a notion of Divergence function, the most notable example of which is commonly known as 
Kullback-Leibler divergence (also called relative entropy) \cite{KL51}. A divergence function is a two-point function 
$$F: \mathcal{S}\times\mathcal{S}\to \mathbb{R}$$
which is differentiable in both entries and with the following properties:
\begin{align}
    &F(\mathbf{x}, \mathbf{y}) \geq 0 . \notag\\
    &F(\mathbf{x}, \mathbf{y}) = 0 \Leftrightarrow \mathbf{x} = \mathbf{y} .\label{eq:propertiesdivergence}\\
    &\frac{\partial F}{\partial x^i }\bigr|_{\mathbf{x} = \mathbf{y}} = \frac{\partial F}{\partial y^j }\bigr|_{\mathbf{x} = \mathbf{y}} = 0, 
\notag
    \end{align}
{for $i,j = 1, \dots, \dim(\mathcal{S})$, and where $\mathbf{x} = \{x^1, x^2, \dots\}$ and  $\mathbf{y} = \{y^1, y^2, \dots\}$ are 
coordinates of the manifold $S$}. The last condition just requires that the first derivatives of the function $F$ vanish on the diagonal. Notice that the three conditions together amount to require that the function $F$ has a local extremum on the diagonal.

In the context of information geometry and of statistical inference the divergence functions are measures of the relative difference 
between probability distributions. We refer to 
\cite{AN00} for an introduction and selected applications of divergence functions in Information Geometry. One remarkable 
property, from the geometric viewpoint, of these general divergence functions is that they are \emph{generating potentials} for the 
geometric structures introduced so far. The concrete meaning of the term ``generating potential'' will be made clearer in the next 
section. The properties of Eq.~\eqref{eq:propertiesdivergence} are sufficient for a two-point function to generate a Riemannian 
metric and a dualistic structure and therefore can be applied in general for any differential manifold, regardless if it is a Statistical 
Manifold or not. The particular case of the Kullback-Leibler divergence gives rise to the Fisher-Rao metric and the so called 
exponential connection and mixture connections.

Having a generating potential for geometric quantities is an advantage, for instance, when deriving geometric or, more properly, 
tensorial quantities that must be preserved under concrete transformations. Indeed, this is the case of the Kullback-Leibler 
divergence, where its monotonicity under transformations fo the probability distributions reflects on the monotonicity of the 
Fisher-Rao metric \cite{Csiszar67, Csiszar67b, Vajda89}.

In the context of Quantum Mechanics, the geometric description of the space of pure states \cite{EMM10, AS99} leads to a natural 
interpretation of the Fubini-Study metric as a quantum counterpart of the Fisher-Rao metric \cite{Wootters81, FN95, Fea10}. 
However, when one wants to apply the same information geometric approach to the space of quantum mixed states the picture is 
not so clear. There are intrinsic difficulties in the problem, like the ambiguity in defining the logarithmic derivative of a matrix 
\cite{Helstrom76, ES13} and the fact that the space of mixed states is not a manifold but it is partitioned into the disjoint union of differential manifolds labelled by the rank of quantum states  \cite{MMZS05, GKM05, Cea17c}. This leads to a wide variety of quantum divergence functions and quantum Fisher-Rao 
metrics depending on the particular aspect that one wants to highlight, cf.\ \cite{Petz96,Cea17,MMVV17}. Having a geometric description of the 
structures involved, independent of the particular realisation or setting, becomes useful in identifying the most relevant among all 
the possibilities.

The aim of this article is to analyse thoroughly the role of two-point functions in differential geometry as potential functions for 
tensorial quantities and we will do so in an intrinsic approach. Our interest relies also in the definition of a well-posed inverse 
problem for the potential function, cf. \cite{Cea17b, CCM17}. In particular we will see that Hamilton's principal function is a particular 
solution to the problem.


\section{Generating Potentials for Tensorial Quantities}

Generating potentials are common objects in mathematics. The most simple intrinsic example that one might think of would be the 
exterior differential of a smooth function $f$ on a differentiable manifold $\M$. The one-form obtained this way is intrinsically 
defined and does not depend on the particular coordinates chosen to define it. At this point one can also recognise that
there might be some global problems when defining the inverse problem. Poincaré's Lemma shows that the inverse problem, i.e., 
if 
given a closed form $\alpha$ one is able to find a function on $\M$ such that $\mathrm{d}f = \alpha$, has a local solution but in 
general a global solution might not exist. 

An interesting situation is encountered in the context of Kähler manifolds where one can describe (at least locally) the symplectic form $\omega$ using the so-called (local) Kähler potential $f$ by means of the formula
\begin{equation}\label{eq:kahler}
    \omega_{ij} = \frac{\partial^2 f}{\partial z\partial\bar{z}},
\end{equation}
where $z$, $\bar{z}$ are, respectively, holomorphic and anti-holomorphic coordinates. Again, this definition is intrinsic and the symplectic form defined this way does not 
depend on the particular coordinates chosen to define it. Of course, by formula \eqref{eq:kahler} one computes the coordinate 
expression of the symplectic form, but this coordinate expression transforms properly under a change of coordinates of the 
manifold. In fact, this situation is analogous to that of the exterior differential. One can show that a local Kähler potential always 
exist. One can easily define the coordinates of a tensor by taking iterated derivatives of a function defined on the manifold. For 
instance, in \cite{Duistermaat01} there is introduced the concept of  \emph{Hessian Riemannian structure}. Let $\mathbf{x} = \{x^1,x^2,\dots\}$ be local coordinates on an open neighbourhood $U\subset\M$ and let $f:U\to\mathbb{R}$ be a convex function. Then one can define the local coordinates of a Riemannian metric by
$$g_{ij}(\mathbf{x}):=\frac{\partial^2 f}{\partial x^i \partial x^j}.$$
While this can be done to define the local coordinates of a tensor, if one chooses to express the function $f$ in a different local chart $\mathbf{y}=\{y^1,y^2,\dots\}$, the expression
$$\tilde{g}_{ij}(\mathbf{y}):=\frac{\partial^2 f}{\partial y^i \partial y^j}$$
defines a different tensor. There is a notable exception, closely related to Kählerian functions, the so-called (flat) \emph{Hessian Structure} \cite{Shima07}. Hence, recovering desired transformation properties of $g$ by imposing invariance conditions on $f$ is not possible by such an approach. Remarkably, if one considers two-point functions with conditions \eqref{eq:propertiesdivergence} one is able to derive tensors of rank higher than one. This is well known for the information geometry community, where the relevant divergence functions that arise give rise to the Fisher-Rao metric. However, it turns out to be a general situation that can be defined on an arbitrary differentiable manifold.

Since we want to define the potential functions in an intrinsic way, we need to define the restriction to the diagonal appropriately. Consider the diagonal embedding given by
\begin{align*}
    D:\M &\hookrightarrow \M\times\M\\
    p &\mapsto (p,p)\qquad.
\end{align*}
The restriction to the diagonal is given by means of the pull-back by $D$. If \hbox{$F:\M\times\M\to\mathbb{R}$} is a two-point function, then $F|_{\mathrm{diag}}= D^*F$. In local coordinates one has that
$$D^*F(p) = F\left(D(p)\right) = F\left(p(x^1,x^2,\dots),p(x^1,x^2,\dots)\right), p\in\M.$$
In addition, one can define the two canonical projections:
\begin{align*}
    \pi_L:\M\times\M &\hookrightarrow \M    &     \pi_R:\M\times\M &\hookrightarrow \M\\
    (p,q) &\mapsto p & (p,q) &\mapsto q\quad.
\end{align*}
By means of these one can define the Left-lift  and Right-lift of a vector field $X\in\mathfrak{X}(\M)$ as follows \cite{Cea17}. The Left-lift $X_L\in\mathfrak{X}(\M\times\M)$ of the vector field $X\in\mathfrak{X}(\M)$ is defined as the unique $\pi_L$-related vector field such that \hbox{$X_L(\pi^*_Rf)=0$} and $X_L(\pi^*_Lf)=\pi^*_L\left(\strut X(f)\right)$. The Right-lift is defined analogously.

We will say that a two-point function 
$$S:\M\times\M\to\mathbb{R}$$
is a potential function if it satisfies the property
\begin{equation}\label{eq:contrast}
D^*(\mathcal{L}_XS) = 0, \quad X\in\mathfrak{X}(\M\times\M)\,,
\end{equation}
where $\mathcal{L}_X$ is the Lie derivative with respect to the vector field $X$.
The condition in \eqref{eq:contrast} is slightly more general than the conditions in \eqref{eq:propertiesdivergence}. Indeed, the latter imply that the diagonal is a local extremum of the function while \eqref{eq:contrast} only requires it to be a critical point.

Any potential function defines a second order tensor in the following way:
$$g(X,Y) := D^*(\mathcal{L}_{X_L}\mathcal{L}_{Y_L}S).$$
We will show that $g$ is $f$-linear in both entries and therefore that it is a tensor. First notice that for any function $f$ on $\M$ and any vector field $X\in\mathfrak{X}(\M)$ one has that
$$(fX)_L = (\pi^*_Lf) X_L$$
 and therefore $f$-linearity in the first entry follows by the $f$-linearity of the Lie derivative. Now consider 
\begin{align}
    \mathcal{L}_{X_L}\mathcal{L}_{(fY)_L}S &= \mathcal{L}_{X_L}\left(\pi^*_Lf\mathcal{L}_{Y_L} S\right)\\
        &= \pi^*_Lf\mathcal{L}_{X_L}\mathcal{L}_{Y_L}S + \mathcal{L}_{X_L}(\pi^*_Lf)\mathcal{L}_{Y_L} S,
\end{align}
and taking the pull-back by $D$ on both sides we get
\begin{equation}
    g(X,fY) = fg(X,Y) + D^*\left( \mathcal{L}_{X_L}(\pi^*_Lf)\right)D^*\left(\mathcal{L}_{Y_L} S\right)
\end{equation}
and by condition \eqref{eq:contrast} the second term at the right hand side vanishes. It is easy to see that this tensor is symmetric. Indeed, 
$$\mathcal{L}_{X_L}\mathcal{L}_{Y_L}S = \mathcal{L}_{Y_L}\mathcal{L}_{X_L}S + \mathcal{L}_{[X,Y]}S$$
and after applying $D^*$ one gets
$$g(X,Y) = g(Y,X).$$

Notice that no assumption has been made on the symmetry of $S$ under permutation of its entries and in the most common situations this will not be the case. There are other possible combinations in order to obtain tensors of rank two, but all of them lead to the same symmetric tensor up to a sign. Indeed,
\begin{align}\label{eq:generateg}
D^*(\mathcal{L}_{X_L}\mathcal{L}_{Y_L}S) &= D^*(\mathcal{L}_{X_R}\mathcal{L}_{Y_R}S) \notag\\
        &= -D^*(\mathcal{L}_{X_R}\mathcal{L}_{Y_L}S) \\
        &= -D^{*}\left(\mathcal{L}_{X_{L}}\,\mathcal{L}_{Y_{R}}\,S\right)=g(X,Y).\notag
\end{align}
A general and coordinate free proof of this fact can be found in \cite{Cea17}.
Note that if every point on the diagonal of $\M\times\M$ is a minimum for $S$, then $g$ is positive semi-definite.
Analogously, if every point on the diagonal of $\M\times\M$ is a maximum for $S$, then $g$ is negative semi-definite.
If $\{x^{j},\,y^{k}\}$ is a coordinate system on $\M\times\M$ adapted to the Cartesian product structure, then the local components of $g$ read:
\be
g_{jk}=\left(\frac{\partial^{2}S}{\partial x^{j}\partial x^{k}}\right)_{x=y} = \left(\frac{\partial^{2}S}{\partial y^{j}\partial y^{k}}\right)_{x=y}= - \left(\frac{\partial^{2}S}{\partial x^{j}\partial y^{k}}\right)_{x=y}\,,
\ee
in complete accordance with the expressions used in information geometry \cite{Aea87, AN00}.

\begin{example}
Let $\mathcal{H}$ be a finite-dimensional complex Hilbert space, and consider the open submanifold $\mathcal{H}_{0}$ consisting of $\mathcal{H}$ without the null vector.
It is well known that the space of pure quantum states may be identified with the complex projective space $\mathbb{CP}(\mathcal{H})$ associated with $\mathcal{H}$.
Furthermore, the complex projective space is the base space of the principal fibre bundle $\pi\colon\mathcal{H}_{0}\rightarrow\mathbb{CP}(\mathcal{H})$ with structure group $\mathbb{C}_{0}$, that is, the Abelian multiplicative group of non-zero complex numbers.
In particular, we have $\pi(\psi)=[\psi]$, where $[\psi]$ is the equivalence class of vectors in $\mathcal{H}_{0}$ differing only by the multiplication with a non-zero complex number.
Now, let us consider the following two-point function on $\mathcal{H}_{0}\times\mathcal{H}_{0}$:
\be
S(\psi,\,\phi):=\frac{|\langle\psi|\phi\rangle|^{2}}{|\psi|^{2}\,|\phi|^{2}}\,,
\ee
which is the quantum-mechanical counterpart of the (non-normalised) ``generalised transition probability'' introduced by Cantoni in \cite{Cant75, Cant77, Wootters81}.
If $\{\mathbf{e}_{j}\}_{1,...,\mathrm{dim}(\mathcal{H})}$ is an orthonormal basis in $\mathcal{H}$, we may introduce a Cartesian coordinate system on $\mathcal{H}_{0}$ setting $\psi=\sum_{j}\,(x^{j} + i y^{j})\,\mathbf{e}_{j}$.
Then, in an obvious way, we may introduce a Cartesian coordinate system $\{x^{j},\,y^{j};\,X^{k},\,Y^{k}\}$ on $\mathcal{H}_{0}\times\mathcal{H}_{0}$ that is adapted to the product structure.
In this (global) coordinate system, it is a matter of straightforward computations to prove that $S$ is a potential function.
Furthermore, from the expressions:
$$
\frac{\partial S}{\partial x^{j}}= \frac{1}{|\psi|^{2}\,|\phi|^{2}}\,\left(\langle\psi|\phi\rangle\,(X^{j} - i Y^{j}) + \langle\phi|\psi\rangle\,(X^{j} + i Y^{j}) - \frac{2x^{j}\,|\langle\psi|\phi\rangle|^{2}}{|\psi|^{2}}\right)\,,
$$
$$
\frac{\partial S}{\partial y^{j}}= \frac{1}{|\psi|^{2}\,|\phi|^{2}}\,\left(\langle\psi|\phi\rangle\,(Y^{j} + i X^{j}) + \langle\phi|\psi\rangle\,(Y^{j} - i X^{j}) - \frac{2y^{j}\,|\langle\psi|\phi\rangle|^{2}}{|\psi|^{2}}\right)\,,
$$
it is easy to see that the symmetric covariant tensor $g$ generated by $S$ reads:

\be
\begin{split}
g =  & \sum_{j,k}\,\frac{2\left(x^{j}x^{k} + y^{j}y^{k} - \delta_{jk}\,R^{2}\right)}{R^{4}}\,\left(\mathrm{d}x^{j}\,\otimes\,\mathrm{d}x^{k} + \mathrm{d}y^{j}\,\otimes\,\mathrm{d}y^{k}\right) +  \\
& +  \sum_{j,k}\, \frac{y^{j}x^{k} - y^{k}x^{j}}{R^{4}}\,\left(\mathrm{d}x^{j}\,\otimes\,\mathrm{d}y^{k} - \,\mathrm{d}y^{j}\,\otimes\,\mathrm{d}x^{k}\right) ,
\end{split}
\ee
where $R^2 = \mathbf{x}^2 + \mathbf{y}^2$. This is the pull-back to $\mathcal{H}_{0}$ of the Fubini-Study metric tensor on the space of pure quantum states $\mathbb{CP}(\mathcal{H})$.
\end{example}

We may take a further step and give a generating algorithm for covariant tensors of order three. The main ingredient in order to do so is that one can combine the action of Left-lifts and Right-lifts in order to cancel non-tensorial contributions of the higher order derivatives. Consider for example
\begin{align}\label{eq:skewness}
    T_1(X,Y,Z) &:= D^*\left(\mathcal{L}_{X_L}\mathcal{L}_{Y_R}\mathcal{L}_{Z_R}S - \mathcal{L}_{X_R}\mathcal{L}_{Y_L}\mathcal{L}_{Z_L}S\right),\\
    T_2(X,Y,Z) &:= D^*\left(\mathcal{L}_{X_L}\mathcal{L}_{Y_L}\mathcal{L}_{Z_R}S - \mathcal{L}_{X_R}\mathcal{L}_{Y_R}\mathcal{L}_{Z_L}S\right).\notag
\end{align}
It is easy to check that both objects are indeed $f$-linear in all the entries and symmetric. One can define up to 8 different possible combinations that lead to order three tensors. As happens for the second order tensor, it can be proven that all of them are equal up to a sign,  cf.\ \cite{Cea17}.

In the context of information geometry this tensors are called the skewness tensors \cite[Chapter4]{Aea87} and they are used to define the dualistic structure of the statistical manifold. The Christoffel symbols of the dual connections $\nabla$ and $\nabla^*$ are defined by adding and subtracting, respectively, to the Christoffel symbol of the Levi-Civita connection the skewness tensor $$\Gamma^{i}_{\negthickspace\pm jk} = \Gamma^{i}_{\negthickspace {jk}} \pm g^{il}T_{ljk}.$$
For instance, by this procedure the Kullback-Leibler divergence gives rise to the exponential connection and the mixture connection. Notice that both right hand sides in \eqref{eq:skewness} vanish if the potential function is symmetric, thus recovering the well known result that symmetric divergence functions can lead only to the self-dual case $\nabla=\nabla^*$.

At this point, a natural question would be if it were possible to use the same procedure to recover higher order tensors from a given potential function $S$. However, before considering this issue in Section \ref{sec:generalcase}, we will address first the inverse problem of finding a potential function that generates two given covariant tensor fields of order two and three.


\section{Hamilton-Jacobi approach to Divergence Functions and the inverse problem}

In this section we will address the inverse problem of finding a potential  function given a Riemannian metric and a skewness tensor. This inverse problem will be relevant in situations where there are not clear candidates for the potential/divergence functions but where there are relevant metrics or connections that one would want to generate. This is of particular importance in Quantum Mechanics where, as explained in the introduction, it can help to distinguish or highlight particular families of divergence functions among the available variety.

In some situations one can find a canonical divergence function. For instance, the self-dual case has as natural divergence function the geodesic distance between points of the manifold. In general, a canonical solution to the problem will not exist. Roughly speaking, this is due to the fact that the metric and skewness tensor only fix the lower order terms in a power series expansion. Thus, the higher order terms might be changed arbitrarily to define infinite possible \emph{equivalent} divergence functions. By equivalent we mean that they generate the same metric and skewness tensors. A successful attempt to obtain solutions to the inverse problem is given, for instance, in \cite{AA15}, however, this approach lacks of a characterisation of the geometric conditions needed to guarantee the global existence of the solution. 

We will review here the approach introduced in \cite{Cea17b} based on the Hamilton-Jacobi theory. In this case the existence of the solution of the inverse problem can be characterised. Moreover, it will allow us to address the problem of finding a generating potential for higher order rank tensors. The main idea is to consider Hamilton's Principal function as a two-point function. Given a Lagrangian $\mathfrak{L}$ and a solution of the equations of motion determined by the Lagrangian, $\gamma:[t_\mathrm{in} , t_{\mathrm{fin}}] \to \M$ Hamilton's principal function is given by the action integral
$$S[\gamma] = \int_{t_\mathrm{in}}^{t_{\mathrm{fin}}} \mathfrak{L}\left( \gamma(t),\dot{\gamma}(t)\right) \mathrm{d}t.$$
If one considers only solutions with $t_\mathrm{in} = 0$ and $t_{\mathrm{fin}} = 1$ such that $\gamma(0)= \mathbf{x}$ and $\gamma(1)= \mathbf{y}$ it becomes a two point function
$$S[\mathbf{x}, \mathbf{y}] = \int_0^1\mathfrak{L}\left( \gamma(t),\dot{\gamma}(t)  \right)\mathrm{d}t.$$
Moreover, assume that the dynamical system defined by the Lagrangian is completely integrable. 

From now on we will assume that the Lagrangian is at least of second order with respect to the velocities, i.e., 
$$\mathfrak{L} = g_{ij}v^iv^j + \mathcal{O}(|v|^3),$$
and that $g$ is positive defined. In particular, let us consider the Lagrangian
\begin{equation}
    \mathfrak{L}(\mathbf{q},\mathbf{v}) = \frac{1}{2} g_{ij}v^i v^j + \frac{\alpha}{6} T_{ijk} v^i v^j v^k,
\end{equation}
where $\alpha\in\mathbb{R}$ is a parameter that is included for convenience. Have into account that this Lagrangian can be written in a coordinate free way as was done in \cite{Cea17b}.

We will show that the two-point function $S$ determined by this Lagrangian solves, under certain assumptions, the inverse problem for $g$ and $T$. 
From the Hamilton-Jacobi theory \cite{Lanczos}, having into account that $S$ is Hamilton's principal function and recalling that in the Lagrangian formulation of dynamics \cite{Mea90} the local expression of the canonical momentum is given by $p_{i}=\frac{\partial \mathfrak{L}}{\partial v^{i}}$, we have that:%
\begin{equation}
    \frac{\partial S}{\partial x^i} = -p_i = \frac{\partial \mathfrak{L}}{\partial v^i}\Bigr|_{\mathbf{q}=\x} = -g_{ij}(\x) v^j_{\i} - \frac{\alpha}{2} T_{ijk}(\x) v^j_\i v^k_\i ,
\end{equation}
\begin{equation}
    \frac{\partial S}{\partial y^i} = P_i = \frac{\partial \mathfrak{L}}{\partial v^i}\Bigr|_{\mathbf{q}=\y} = g_{ij}(\y) v^j_{\f} + \frac{\alpha}{2} T_{ijk}(\y) v^j_\f v^k_\f .
\end{equation}
This guarantees that 
$$D^*p_i = D^* P_i =0,\quad i= 1,\dots,n.$$
This is because the points in the diagonal satisfy $\x=\gamma(0)=\gamma(1)=\y$ and therefore the initial velocity is zero. This shows that the two-point function satisfies the condition in Eq.~\eqref{eq:contrast}.

We need to obtain the explicit solution of the equations of motion in order to be able to express $\mathbf{v}_\i$ and $\mathbf{v}_\f$ only in terms of the endpoints of the curve $\gamma$. If we use the Taylor expansion around $t=0$ of the solution of the equations of motion we have that
\begin{equation}
    \gamma(t) = \gamma(0) + t\mathbf{v}_\i + \frac{t^2}{2}\dot{\mathbf{v}}_\i + \frac{t^3}{6}\ddot{\mathbf{v}}_\i + \mathcal{O}(t^4).
\end{equation}
Substituting $t=1$ and having into account that $\gamma(0) = \x$ and $\gamma(1) = \y$ we have that
\begin{equation}\label{eq:vin}
    \mathbf{v}_\i = \y - \x - \frac{1}{2}\dot{\mathbf{v}}_\i - \frac{1}{6} \ddot{\mathbf{v}}_\i + \mathcal{O}(t^4).
\end{equation}
Analogously, if we expand around $t=1$ we get an expression for $\mathbf{v}_\mathrm{fin}$, namely,
\begin{equation}\label{eq:vfin}
    \mathbf{v}_\f = \y - \x + \frac{1}{2}\dot{\mathbf{v}}_\f - \frac{1}{6} \ddot{\mathbf{v}}_\f + \mathcal{O}(t^4).
\end{equation}
The Euler-Lagrange equations determined by the Lagrangian are
\begin{equation}\label{eq:eulerlagrange}
    \dot{v}^l = -\alpha T^{l}_{jk}v^j\dot{v}^k - \Gamma^{l}_{jk} v^j v^k - \frac{\alpha}{6} g^{lr} A_{rjks}v^j v^k v^s,
\end{equation}
where we are raising and lowering indices in the usual way by means of the Riemannian metric $g$. Here, $\Gamma^{i}_{jk}$ are the Christoffel symbols of the Levi-Civita connection and we have defined the tensor
$$A_{rjks}:= \left( \frac{\partial T_{jkr}}{\partial q^{s}} + \frac{\partial T_{jrs}}{\partial q^{k}} + \frac{\partial T_{rks}}{\partial q^{j}} - \frac{\partial T_{jks}}{\partial q^{r}} \right).$$

We will assume now that $\dot{\mathbf{v}}$ is an analytic function of the velocities and we will neglect terms with high powers of $\mathbf{v}$, since we are interested in what happens in a neighbourhood of $\mathbf{v} = 0$, i.e.\ the diagonal. We can inductively substitute \eqref{eq:eulerlagrange} in \eqref{eq:vin} and \eqref{eq:vfin} in order to obtain the expressions for $\mathbf{v}_\i$ and $\mathbf{v}_\f$ respectively. The analytic expression for $\ddot{\mathbf{v}}$ is obtained by differentiating the equations of motion. Notice that $\dot{\mathbf{v}}$ is at least of order $\mathcal{O}(\mathbf{v}^2)$ and that $\ddot{\mathbf{v}}$ is at least of order $\mathcal{O}(\mathbf{v}^3)$. If we define now $\Delta^i = y^i - x ^i$ a straightforward calculation leads to:
\begin{align}
    \frac{\partial S}{\partial x^i} = &-g_{ij}(\x) \Delta^j -\frac{1}{2}\Gamma_{ijk}(\x) \Delta^j \Delta^k -\frac{1}{6} \Gamma_{ijk}(\x) \Gamma^{j}_{ls}(\x) \Delta^k \Delta^l \Delta^s \label{eq:1stderivativea}\\
        &-\frac{\alpha}{2} T_{ijk}(\x) \Delta^j \Delta^k -\frac{\alpha}{12} A_{ijkl}(\x) \Delta^j \Delta^k \Delta^l + \mathcal{O}\left( (\y-\x)^4\right).\notag
\end{align}
\begin{align}
    \frac{\partial S}{\partial y^i} = &\phantom{-}g_{ij}(\y) \Delta^j -\frac{1}{2}\Gamma_{ijk}(\y) \Delta^j \Delta^k +\frac{1}{6} \Gamma_{ijk}(\y) \Gamma^{j}_{ls}(\y) \Delta^k \Delta^l \Delta^s\label{eq:1stderivativeb}\\
        &+\frac{\alpha}{2} T_{ijk}(\y) \Delta^j \Delta^k -\frac{\alpha}{12} A_{ijkl}(\y) \Delta^j \Delta^k \Delta^l + \mathcal{O}\left( (\y-\x)^4\right).\notag
\end{align}
Now one can check that formulae \eqref{eq:generateg} and \eqref{eq:skewness} give rise to the tensors $g$ and $T$. For instance, one has that
$$D^*\left( \frac{\partial ^2 S}{ \partial y^i \partial x^j} \right)(\mathbf{q}) = D^*\left( \frac{\partial ^2 S}{ \partial x^i \partial y^j} \right)(\mathbf{q}) = -g_{ij} = - D^*\left( \frac{\partial ^2 S}{ \partial x^i \partial x^j} \right)(\mathbf{q}).$$
this shows that (up to a sign) only one tensor of rank two can be generated out of the two-point function. Moreover,
\begin{equation}\label{eq:3rdderivativesa}
    D^*\left( \frac{\partial ^3 S}{  \partial x^l \partial x^k \partial y^j} \right)(\mathbf{q}) = -\Gamma_{jkl} + \alpha T_{jkl},
\end{equation}
\begin{equation}\label{eq:3rdderivativesb}
    D^*\left( \frac{\partial ^3 S}{  \partial y^l \partial y^k \partial x^j} \right)(\mathbf{q}) = -\Gamma_{jkl} - \alpha T_{jkl}.
\end{equation}
Hence, choosing $\alpha = 1/2$ in the Lagrangian solves the inverse problem, since subtracting the latter formulae from each other is precisely one of the definitions in \eqref{eq:skewness}. Remember that there are other 6 possible linear combinations of the derivatives that lead to rank three tensors. However, one can check that they generate the same tensor $T$ up to a sign.

Strictly speaking, in order to prove \eqref{eq:3rdderivativesa} and \eqref{eq:3rdderivativesb} it is not necessary to consider the term $\ddot{\mathbf{v}}$ in the expressions \eqref{eq:vin} and \eqref{eq:vfin}. Neither is it necessary to consider the expressions \eqref{eq:1stderivativea} and \eqref{eq:1stderivativeb} up to order $\mathcal{O}((\y-\x)^3)$. The contributions of these terms to the third order derivatives of the principal function vanish identically after taking the pull-back to the diagonal. In fact, these terms were not taken into account in \cite{Cea17b} but have been considered here since we are willing to study what happens in order to generate higher order tensors. In particular, we are interested in the possibility of generating fourth order tensor fields.


\section{Higher rank tensors in the general case}\label{sec:generalcase}

In this section we will consider the problem of defining tensors of order four out of a two-point function, or more concretely, out of a potential function, cf.\ Eq.~\eqref{eq:contrast}. Our first objective will be to find linear combinations of derivatives of the two-point function giving raise to rank four tensors. We will use the following abbreviated notation
$$LRLR := D^* \left(\mathcal{L}_{X_L}\mathcal{L}_{Y_R}\mathcal{L}_{Z_L}\mathcal{L}_{W_R}S\right).$$
And similarly for other combinations of Left-\ and Right-lifts. This is convenient because the only differences between the terms will be the particular combination of Left-\ and Right-lifts that appear. There are only two possible rank four tensors, namely
\begin{align}
Q_1(X,Y,Z,W) = &(LRRL - RLLR) + (LRRR - RLLL)\notag\\
        &+ (LRLL - RLRR) + (LRLR - RLRL),
\end{align}
\begin{align}
Q_2(X,Y,Z,W) = &(LLLL-RRRR) + (LLLR - RRRL)\notag\\
        &+ (LLRL - RRLR) + (LLRR - RRLL).
\end{align}
We can now compute the derivatives of the principal function and see what is the rank four tensor that one gets. In principle, one could also enlarge the Lagrangian to accommodate a term of order $\mathcal{O}(\mathbf{v}^4)$. Before computing the tensors $Q_1$ and $Q_2$ we will consider the Lagrangian given by
\begin{equation}
    \mathfrak{L}(\mathbf{q},\mathbf{v}) = \frac{1}{2} g_{ij}v^i v^j + \frac{\alpha}{6} T_{ijk} v^i v^j v^k +\frac{1}{24} C_{ijkl}  v^i v^j v^k v^l,
\end{equation}
where now $C$ is some fourth order tensor. Using the same analysis of the previous section one can show that the first order derivatives of the principal function become now:
\begin{align}
    \frac{\partial S}{\partial x^i} = &-g_{ij}(\x) \Delta^j -\frac{1}{2}\Gamma_{ijk}(\x) \Delta^j \Delta^k -\frac{1}{6} \Gamma_{ijk}(\x) \Gamma^{j}_{ls}(\x) \Delta^k \Delta^l \Delta^s \label{eq:1stderivativea2}\\
        &-\frac{\alpha}{2} T_{ijk}(\x) \Delta^j \Delta^k -\frac{\alpha}{12} A_{ijkl}(\x) \Delta^j \Delta^k \Delta^l \notag\\
        &-\frac{1}{24} C_{ijkl}(\x)  \Delta^j \Delta^k \Delta^l+ \mathcal{O}\left( (\y-\x)^4\right).\notag
\end{align}
\begin{align}
    \frac{\partial S}{\partial y^i} = &\phantom{-}g_{ij}(\y) \Delta^j -\frac{1}{2}\Gamma_{ijk}(\y) \Delta^j \Delta^k +\frac{1}{6} \Gamma_{ijk}(\y) \Gamma^{j}_{ls}(\y) \Delta^k \Delta^l \Delta^s\label{eq:1stderivativeb2}\\
        &+\frac{\alpha}{2} T_{ijk}(\y) \Delta^j \Delta^k -\frac{\alpha}{12} A_{ijkl}(\y) \Delta^j \Delta^k \Delta^l\notag\\
        &+\frac{1}{24} C_{ijkl}(\y)  \Delta^j \Delta^k \Delta^l+ \mathcal{O}\left( (\y-\x)^4\right).\notag
\end{align}
Now one can compute the result of tensors $Q_1$ and $Q_2$ and it turns out that both of them vanish identically. This means that this particular potential function does not generate rank four tensors.

In what follows we will show that if the potential function is analytic on a neighbourhood of the diagonal, it cannot generate order four tensors. First notice that condition \eqref{eq:contrast} amounts to say that the diagonal of $\M\times\M$ is a set of critical points for a potential function $S$, and there it attains the value $S(\mathbf{q},\mathbf{q})$ for $\mathbf{q}\in\M$. A Taylor series expansion around a point $\mathbf{q}_0$ of the diagonal has necessarily the form
\begin{align}
    S(\x,\y) = & S(\mathbf{q_0},\mathbf{q_0}) + h_{ij}(\mathbf{q}_0)(y^i - x^i)(y^j - x^j) + \notag \\
    & + t_{ijk}(\mathbf{q}_0)_{ijk} (y^i - x^i)(y^j - x^j) (y^k - x^k) + \\
            &+ c_{ijkl}(\mathbf{q}_0) (y^i - x^i)(y^j - x^j) (y^k - x^k)(y^l - x^l) + \mathcal{O}((\y - \x)^5),\notag
\end{align}
where $\{x^1, x^2, \dots, x^n, y^1, y^2,\dots, y^n \}$ are now coordinates relative to the point $(\mathbf{q}_0,\mathbf{q}_0)$ on $\M\times\M$ and where $h_{ij}$ depends on the second derivatives of $S$, $t_{ijk}$ depends on the third derivatives of $S$ and so on. This is so because all the derivatives along the direction tangent to the diagonal vanish.

If one computes the tensors $Q_1$ and $Q_2$ in this particular case one gets that they vanish identically. This means that a potential function, assuming it is analytic, is not capable of generating non-trivial covariant tensor fields of order four.


\section{Conclusions}

We have studied the inverse problem for potential functions and found a solution in terms of Hamilton's principal function and a suitable chosen Lagrangian built from the metric and the skewness tensor. The solution will exist  provided that the equations of motion of the associate Hamilton-Jacobi problem are completely integrable. We have also studied if one is able to recover higher order rank tensors (order four) from a potential function satisfying \eqref{eq:contrast}. It turns out that the answer is negative assuming that the potential fucntion is analytic on a neighbourhood of the diagonal.

Notice that potential functions, as they are described here, have the limitation that they cannot generate non-symmetric tensors. For instance, one cannot generate a symplectic structure using the procedure described. One way of achieving that could be by considering Grassmann variables. 

Another interesting line to pursue would be that of finding a way of generating higher order tensors. From what we have seen, at least for order four, this will not be achievable by a two-point function. By that we mean that the tensor arises as linear combinations of the derivatives of the potential function. It is clear that once one has obtained the Christoffel symbols, for instance, one can compute out of them the Riemann tensor in the usual way. What we have shown is that the Riemann tensor cannot be directly obtained by computing fourth order derivatives of the potential function. It is very likely that one could do that by considering three-point functions. However, developing such a theory in depth would also require a Hamilton-Jacobi theory of order higher than two (associated to evolution equations including up to third order derivatives).



\section*{Acknowledgements}
 
The authors want to thank the scientific committee of the ``XXVI International Fall Workshop on Geometry and Physics'' for giving them the opportunity to present this research. G.M. acknowledges financial support from the Spanish Ministry of Economy and Competitiveness, through the Severo Ochoa Programme for Centres of Excellence in RD (SEV-2015/0554), and would like to thank the support provided by the Santander/UC3M Excellence Chair Programme 2016/2017. J.M.P.P. was supported by QUITEMAD+ S2013/ICE-2801, the Spanish MINECO grant MTM2017-84098-P and the ``Juan de la Cierva - Incorporaci\'on" Proyect 2018/00002/001.


\end{document}